\documentclass[preprint,12pt]{elsarticle}
\usepackage{amssymb, amsmath, amsthm}
\usepackage{color}

\biboptions{sort&compress}

\newtheorem{theorem}{Theorem}[section]

\newtheorem{proposition}[theorem]{Proposition}
\newtheorem{corollary}[theorem]{Corollary}

\newtheorem{definition}[theorem]{Definition}

\numberwithin{equation}{section}

\newcommand\NN{\mathbb{{N}}}

\newcommand\RR{\mathbb{{R}}}
\newcommand\CC{\mathbb{{C}}}
\newcommand\ZZ{\mathbb{{Z}}}

\newcommand\rE{{\rm E}}

\newcommand\ra{{\rm a}}

\newcommand\rd{{\rm d}}

\newcommand\be{{\boldsymbol e}}

\newcommand\bj{{\bf j}}

\newcommand\bell{{\boldsymbol\ell}}

\newcommand\bt{{\bf t}}

\newcommand\bx{{\boldsymbol x}}

\newcommand\bz{{\boldsymbol z}}

\newcommand\balpha{{\boldsymbol\alpha}}
\newcommand\bbeta{{\boldsymbol\beta}}

\newcommand\bveps{{\boldsymbol\epsilon}}
\newcommand\beps{{\boldsymbol\varepsilon}}

\newcommand\btau{{\boldsymbol\tau}}

\newcommand\bnull{{\boldsymbol0}}
\newcommand\bone{{\boldsymbol1}}

\newcommand\cS{{\mathcal S}}

\newcommand{\mfrac}[2]%
{\raisebox{0.5pt}{\footnotesize$\dfrac{#1}{#2}$}}
\newcommand{\mbinom}[2]%
{\raisebox{0.5pt}{\footnotesize$\dbinom{#1}{#2}$}}
\def\smmat\{#1&#2\cr#3&#4\}%
   {\!\!\pmatrix{\!#1\!&\!#2\!\cr\!#3\!&\!#4\!}\!\!}
\newcommand\scrm{{\raise0.5pt\hbox{-}}}
\def\eop{{ \vrule height7pt width7pt depth0pt}\par\bigskip} 
\newcommand\ie{{\it\thinspace i.e.}}
\newenvironment{pf}%
{\par\noindent\textbf{Proof:}~}

\begin{document}

\title{Polynomial Reproduction of Multivariate Scalar Subdivision Schemes\\ with General Dilation}
\author[ch]{Maria Charina\corref{cor1}}
\ead{maria.charina@uni-dortmund.de}
\author[cc]{Lucia Romani}
\ead{lucia.romani@unimib.it}

\cortext[cor1]{Corresponding author}
\address[ch]{Fakult\"at f\"ur Mathematik, TU Dortmund, D--44221 Dortmund, Germany}
\address[cc]{Department of Mathematics and Applications, University of Milano-Bicocca, Via R.Cozzi 53, 20125 Milano, Italy}

\begin{abstract}
In this paper we study scalar
multivariate subdivision schemes with general integer expanding dilation matrix.
Our main result yields simple algebraic conditions on the symbols of such
schemes that characterize their polynomial
reproduction, i.e. their capability to generate
exactly the same polynomials from which the initial data is
sampled. These algebraic conditions also allow us to determine the 
approximation order of the associated refinable functions
and to choose the ``correct" parametrization,
i.e. the  grid points to which the newly computed values are
attached at each subdivision iteration.  We use this special choice of the
parametrization to increase the degree of polynomial reproduction
of known subdivision schemes and  to construct new schemes with given
degree of polynomial reproduction. 
\end{abstract}


\maketitle

\section{Introduction}

In \cite{CharinaConti2012}, the authors derive simple algebraic
conditions on the subdivision symbol of a multivariate scalar
subdivision scheme with dilation matrix $M=mI$, $m \in \ZZ$, $|m| > 1$, that
allow us to determine the degree of polynomial reproduction of
such a scheme and, consequently by \cite{ALevin2003}, its approximation order.
In this paper we extend the results of \cite{CharinaConti2012} to the case of general
 expanding dilation matrix $M \in \ZZ^{s \times s}$ whose spectral
radius $\rho(M)$ satisfies $\rho(M)>1$. On the one hand,
our interest in the study of such a general dilation matrix $M$ is motivated
by the existence of multivariate scalar subdivision schemes for which the results of \cite{CharinaConti2012} are not
applicable, e.g.
$$
 M=\left( \begin{array}{cc} a&b\\0&c\end{array}\right), \quad a, b, c  \in \ZZ, \quad |a|,|c|> 1,
$$
see \cite{Dahlke, GH}.
On the other hand, even for schemes whose dilation matrix $M$ satisfies $M^\ell=mI$ for some integer $\ell>1$
and $m \in \ZZ$, $|m|>1$, e.g., $\sqrt{3}-$subdivision with $M^2=-3I$,
the results of this paper allow us to directly investigate the properties of the associated symbol,
instead of
dealing with its $\ell-$th iterated version.\\
Moreover, in contrast to polynomial generation, i.e. the capability of a scheme to
generate the full space of polynomials of certain degree, there
are very few theoretical results on polynomial reproduction of subdivision schemes,
(see \cite{ContiHormann, ContiRomani, DynShen2008}). Thus, we strongly
believe that our paper may further contribute to the development of this topic.
See \cite{CharinaConti2012} for a detailed
discussion on polynomial generation and reproduction of subdivision schemes.

In addition, from a practical point of view, the results in this paper, if combined with those
in \cite{CharinaContiJetterZimm2010}, provide a simple algebraic tool for the construction of
new interesting multivariate subdivision schemes with enhanced properties. Note that, in the univariate
case, affine combinations of existing schemes have been already used in \cite{CGR2009, CGR2011, ContiRomani2010, ContiRomani2010b, DD, DongShen,
DynShen2008} for such purpose.

\smallskip \noindent  The remainder of this paper is organized as follows. In section \ref{sec:background}, we
briefly summarize key notions concerning subdivision schemes with integer, expanding dilation matrix $M$.
In section \ref{sec:algebra}, we prove that for a non-singular multivariate scalar subdivision scheme with finitely supported mask
$\ra = \{\ra_\balpha,\ \balpha\in\ZZ^s\}$ and symbol
$\displaystyle{a(\bz)=\sum_{\balpha\in\ZZ^s}
\ra_\balpha\;\bz^\balpha}$, the polynomial reproduction of degree up to
$k$ is equivalent to
\begin{equation}\label{eq:main}
 \bigl(D^\bj a\bigr)(1,\dots,1)=\left| \hbox{det}(M)\right| \prod_{i=1}^s \prod_{\ell_i=0}^{j_i-1}(\tau_i-\ell_i)
 \quad \hbox{and} \quad \bigl(D^\bj a\bigr)(\beps)=0 \ \
\end{equation} 
for $\beps\in \Xi'$, $\bj\in \NN_0^s$ and $|\bj| \le k$.
The set $\Xi'$  in \eqref{eq:main} is a finite set of certain multi-indices and
$\btau=(\tau_1,\cdots,\tau_s)\in \RR^s$ appears in the
parametrization associated with the subdivision scheme. The
importance of condition \eqref{eq:main} for $k=1$ is that it
allows us to identify the correct parametrization of any
non-singular or even only convergent subdivision scheme to
guarantee at least the reproduction of linear polynomials. The
parametrization determines the grid points to which the newly
computed values are attached at each step of the subdivision recursion
to ensure the highest degree of polynomial reproduction of the scheme.
Since some preliminary results discussed in \cite{CharinaConti2012} are still valid
for a general integer expanding matrix $M$,
in the following we only prove the extensions of \cite[Propositions 2.3 and 2.5]{CharinaConti2012}.
In section \ref{sec:examples}, we use the results in section \ref{sec:algebra}
to increase the degree of polynomial reproduction of existing subdivision schemes by considering their
affine combinations.

\section{Background and notation}\label{sec:background}

Let $M \in \ZZ^{s \times s}$ be an expanding \emph{dilation matrix}  whose
 spectral radius $\rho(M)$ satisfies $\rho(M)>1$.
The $r-$th step of a scalar $s$-variate subdivision scheme
is given by
\begin{equation}\label{eq:subrecursion}
  \rd^{(r+1)}_\balpha=(\cS_\ra \rd^{(r)})_\balpha=\sum_{\bbeta\in\ZZ^s}
 \ra_{\balpha-M\bbeta}  \rd^{(r)}_\bbeta\;, \quad \rd^{(0)} \in \ell(\ZZ^s), \quad r \ge 0\;, \quad
 \balpha\in\ZZ^s,
\end{equation}
and is defined by the
finitely supported \emph{mask} $\ra = \{\ra_\balpha\in \RR,\
\balpha\in\ZZ^s\}$ and the \emph{subdivision operator}
\begin{equation}\label{eq:subdop}
 \cS_\ra: \ell(\ZZ^s) \rightarrow \ell(\ZZ^s), \quad  \big(\cS_\ra \rd\big)_\balpha = \sum_{\bbeta\in\ZZ^s}
  \ra_{\balpha-M\bbeta}\; \rd_\bbeta\;,\quad \balpha\in\ZZ^s\;,
  \end{equation}
where $\ell(\ZZ^s)$ is the space of scalar sequences indexed by
$\ZZ^s$.
We denote the \emph{symbol} of a scalar $s$-variate subdivision scheme with mask $\ra$ by
\begin{equation*}\label{eq:laurpol}
  a(\bz) = \sum_{\balpha\in\ZZ^s} \ra_\balpha\;\bz^\balpha, \quad \bz = (z_1,\dots,z_s) \in\left( \CC\setminus\{0\}\right)^s,
  \quad \bz^\balpha = z_1^{\alpha_1}z_2^{\alpha_2}\cdot \ldots \cdot z_s^{\alpha_s}.
\end{equation*}
Let $m:=|\hbox{det}(M)|$ and
\begin{equation}\label{eq:extrpts}
  \rE = \{\be_0,\dots,\be_{m-1}\}
\end{equation}
be the set of representatives of $\ZZ^s/ M  \ZZ^s$ containing $\bnull
= (0,0,\dots,0)$. The $m$ \emph{submasks} of $\ra$ and the associated symbols $a_{\be}(\bz)$, $\be \in E$,
are defined by $\{ \ra_{M \balpha+\be}, \, \balpha \in \ZZ^s \}$ and
 $\displaystyle a_{\be}(\bz)=\sum_{\balpha \in \ZZ^s} \ra_{M \balpha+\be}
\, \bz^{M\balpha+\be}$, respectively.

\begin{definition}[Condition Z$_k$]
We say that the mask symbol $a(\bz)$ satisfies the \emph{zero condition of
order $k$ (Condition~Z$_k$)}, if
\begin{equation} \label{eq:Zkcond}
a(\bone)=m \quad \hbox{and} \quad \bigl(D^\bj a\bigr)(\beps)=0
\quad \hbox{for} \quad \beps\in \Xi':=\Xi\setminus \{\bone\} \quad
\hbox{and} \quad |\bj|<k\;,
\end{equation}
where $ \Xi =\{\varepsilon_0, \dots, \varepsilon_{m-1}\}$ is the
set of representatives of  $\ZZ^s / M^T \ZZ^s$ which
contains $\bone=(1, 1, \dots, 1)$.
\end{definition}

For $\btau \in \RR^s$ and $\balpha \in \ZZ^s$ define
\begin{equation}\label{def:general_par}
   \bt_\balpha^{(r)}:= \bt_\bnull^{(r)}+ M^{-r}\balpha,
   \quad \bt_\bnull^{(r)}=\bt_\bnull^{(r-1)}-M^{-r} \btau, \quad \bt_\bnull^{(0)}=0, \quad r\ge 0.
\end{equation}
We call the sequence $\{\bt^{(r)},\ r\ge
 0\}$, $\bt^{(r)}=\{\bt^{(r)}_\balpha,\ \alpha\in \ZZ^s\}$,
the \emph{sequence of parameter values associated with the
subdivision scheme}.

\begin{definition}
If the sequence of continuous functions $\{F^{(r)},\ r\ge 0\}$
with $F^{(r)}( \bt_\balpha^{(r)}) = \rd_\balpha^{(r)}$, $\balpha
\in \ZZ^s$, converges for all initial data
$\rd^{(0)}=\{\rd^{(0)}_\balpha,\ \balpha\in \ZZ^s\} \in
\ell(\ZZ^s)$, then we denote its limit by $
\displaystyle
 S_\ra^\infty \rd^{(0)} = \lim_{r \to\infty} F^{(r)}$ and say that $S_\ra$ is convergent.
\end{definition}

\begin{definition}
A subdivision scheme is called \emph{non-singular}, if it is
convergent, and $S_\ra^\infty \rd^{(0)}=0$ if and only if
$\rd^{(0)}_\balpha=0$ for all $\balpha\in \ZZ^s$.
\end{definition}

We denote by $\Pi_k$ the space of multivariate  polynomials of
total degree $k \ge 0$.

\begin{definition}\label{def:pr}
A convergent subdivision scheme $S_{{\ra}}$ with parametrization $\{\bt^{(r)},\ r\ge 0\}$ is
\emph{reproducing} polynomials up to degree $\gamma_R$ if for any $\pi
\in \Pi_{\gamma_R}$ and  ${\rd}^{(0)}=\{\pi(\bt^{(0)}_\balpha),\
\balpha\in\ZZ^s\}$ the limit of the subdivision process satisfies
$S^{\infty}_{\ra}{\rd}^{(0)}=\pi$.
\end{definition}

\section{Algebraic conditions for polynomial  reproduction} \label{sec:algebra}

In this section we extend the results of Proposition 2.3 and Proposition 2.5 in \cite{CharinaConti2012}
to the case of a general expanding dilation matrix $M \in \ZZ^{s \times s}$.

\begin{proposition} \label{prop:linear-reproduction}
Let $S_\ra$ be a non-singular $s$-variate scalar subdivision scheme that generates
linear polynomials, i.e. its symbol satisfies Condition $Z_2$.
Then $S_\ra$ reproduces linear polynomials if and only if its
parameter values are given by (\ref{def:general_par}) with $$
\btau=m^{-1} \left(D^{\bveps_1} a(\bone),\dots, D^{\bveps_s}
a(\bone)\right),$$
where $D^\bj$ is the $\bj$-th directional derivative and $\bveps_\ell$ is the $\ell$-th
unit vector of $\RR^s$.
\end{proposition}

\pf By \cite[Proposition 1.7]{CharinaConti2012} it suffices to prove
the claim for the step-wise polynomial reproduction. Moreover, any
convergent subdivision scheme reproduces constants, hence we only
consider the starting sequences
$\rd^{(r)}_\balpha=\pi(\bt^{(r)}_\balpha)=(\bt^{(r)}_\balpha)_j$, $\balpha \in \ZZ^s$
and  $j=1, \dots,s$. Then  for any
$\balpha\in\ZZ^s$ and $\be=(\be_1, \dots, \be_s) \in E$ we get
\begin{align*}
  \rd^{(r+1)}_{M\balpha+\be}
  &= \sum_{\bbeta\in\ZZ^s} \ra_{M(\balpha-\bbeta)+\be} \, \rd^{(r)}_\bbeta
   = \sum_{\bbeta\in\ZZ^s} \ra_{M\bbeta+\be} \, \rd^{(r)}_{\balpha-\bbeta} \\
   &= \sum_{\bbeta\in\ZZ^s} \ra_{M\bbeta+\be} \, 
   \left(\bt_{\bnull}^{(r)}+M^{-r}(\balpha-\bbeta)\right)_j \\
  &=  \sum_{\bbeta\in\ZZ^s} \, \ra_{M\bbeta+\be}  \left(\bt_{\bnull}^{(r)}+M^{-(r+1)}(M
  \balpha+\be)\right)_j
     - \sum_{\bbeta\in\ZZ^s} \ra_{M\bbeta+\be} \left(M^{-(r+1)}(M \bbeta+\be)\right)_j\\
     &= \left(\bt_{\bnull}^{(r)}+M^{-(r+1)}(M \balpha+\be)\right)_j
     - \left( M^{-(r+1)} \sum_{\bbeta\in\ZZ^s} \ra_{M\bbeta+\be} (M\bbeta+\be)\right)_j\\
  &= \left(\bt_{\bnull}^{(r)}+M^{-(r+1)}(M \balpha+\be)\right)_j
- \left( M^{-(r+1)}( D^{\bveps_1} a_\be(\bone), \dots,
D^{\bveps_s}
  a_\be(\bone))\right)_j\\
   &= \left(\bt_{\bnull}^{(r)}+M^{-(r+1)}(M \balpha+\be)- m^{-1}
   M^{-(r+1)} ( D^{\bveps_1} a(\bone), \dots, D^{\bveps_s} a(\bone)) \right)_j,
\end{align*}
where we use the fact that $\displaystyle \sum_{\bbeta\in\ZZ^s} \, \ra_{M\bbeta+\be}=a_\be(\bone)=1$ and the last equality is due to
\cite[Proposition 2.1, part $(ii)$ ]{CharinaConti2012}
for $\bj=\bveps_j$. Thus, $
\rd^{(r+1)}_{M\balpha+\be}$ is equal to
\[
  \pi(\bt^{(r+1)}_{M\balpha+\be})=t_{\bnull,j}^{(r+1)}+ \left(M^{-(r+1)} (M \balpha+\be) \right)_j=
  t_{\bnull,j}^{(r)}+ \left(M^{-(r+1)} (M \balpha-\btau +\be) \right)_j,
\]
$\alpha \in \ZZ^s$, 
if and only if $\tau_j=m^{-1} D^{\bveps_j} a(\bone)$  for all $j=1,\cdots,s$. \eop

As in \cite{CharinaConti2012} we easily get the following consequence of
Proposition \ref{prop:linear-reproduction} for convergent subdivision schemes.

\begin{corollary} \label{cor:linear-reproduction}
Let $S_\ra$ be a  convergent $s$-variate scalar subdivision scheme that generates
linear polynomials, i.e. its symbol satisfies Condition $Z_2$.
Then $S_\ra$ reproduces linear polynomials if its parameter values
are given by (\ref{def:general_par}) with
$$
\btau=m^{-1} \left(D^{\bveps_1} a(\bone),\dots, D^{\bveps_s}
a(\bone)\right).
$$
\end{corollary}

The next Proposition is crucial for the proof  of our main result.

\begin{proposition} \label{prop:equiv3}
Let $k \in \NN$, $\btau\in \RR^s$ and $q_{\bj}$ given by
 \begin{equation}\label{def:qj}
q_{\bnull}(\bz):=1,\quad  q_{\bj}(z_1,\dots, z_s):=\prod_{i=1}^s
\prod_{\ell_i=0}^{j_i-1}(z_i-\ell_i), \quad \bj=(j_1,\dots, j_s), \quad \bz \in \RR^s.
 \end{equation}
A subdivision symbol $a(\bz)$ satisfies
\begin{equation}\label{eq:equiv1a}
 \bigl(D^\bj a\bigr)(\bone)=m \  q_{\bj}(\btau), \quad \bigl(D^\bj a\bigr)(\beps)=0
 \quad    \beps\in \Xi',\quad \bj \in \NN_0^s, \quad |\bj| \le k,
\end{equation}
if and only if
\begin{equation}\label{eq:equiv2a}
\sum_{\bbeta\in \ZZ^s} \ra_{\balpha-M\bbeta}
\left(M^{-r}\bbeta\right)^\bj =\left( M^{-(r+1)}(
\balpha-\btau)\right)^\bj, \  \alpha\in \ZZ^s, \  \bj \in
\NN_0^s, \  |\bj|\le k, \  r\ge 0\,.
\end{equation}
\end{proposition}
\pf  Note that due to \cite[Proposition 2.1- part $(iii)$ ]{CharinaConti2012} conditions
in \eqref{eq:equiv1a} are equivalent to
\begin{equation}\label{eq:propqj}
  q_{\bj}(\btau)=\sum_{\bbeta\in \ZZ^s} q_{\bj}(\balpha-M\bbeta)\ra_{\balpha-M\bbeta}, \quad
\bj \in \NN_0^s, \quad |\bj| \le k, \quad \balpha\in\ZZ^s\,.
\end{equation}
The proof is by induction on $k$. For $k=0$ we get, for any $\btau
\in \RR^s$, $ \displaystyle
  q_{\bnull}(\btau)=\sum_{\bbeta\in \ZZ^s}\ra_{\balpha-M\bbeta}=\left(M^{-(r+1)}(
\balpha-\btau)\right)^{\bnull}=1$. Assume next that the claim is
true for all   $\bj\in \NN_0^s$ with  $|\bj| \le k-1$ and prove it for $\bj \in \NN_0^s$
with $|\bj|=k$. The polynomial $q_{\bj}$ in $\bx$ of (total)
degree $|\bj|=k$ is of the form
\begin{equation}\label{eq:propqj2}
q_{\bj}(\balpha-M^{r+1} \bx)=\sum_{\bell \in \NN_0^s, \ |\bell|\le
k} c_{\bj,\balpha,\bell} \bx^\bell,
\quad \bx \in \RR^s\,.
\end{equation}
 Therefore, using the induction assumption  and by \eqref{eq:propqj} and \eqref{eq:propqj2} we have

$$\begin{array}{ll}
    q_{\bj}(\btau)&= \displaystyle \sum_{\bbeta\in \ZZ^s} q_\bj (\balpha-M^{r+1}M^{-r}\bbeta) \ra_{\balpha-M\bbeta}=
     \displaystyle{\sum_{\bbeta\in \ZZ^s}\sum_{\bell  \in \NN_0^s, \ |\bell|\le k}
    c_{\bj,\balpha,\bell} \   \ra_{\balpha-M\bbeta}} \left( M^{-r} \bbeta \right)^\bell \\\\
    &=\displaystyle{\sum_{\bell  \in \NN_0^s, \ |\bell|= k}c_{\bj,\balpha,\bell} \sum_{\bbeta\in \ZZ^s}
    \ra_{\balpha-M\bbeta} \left( M^{-r} \bbeta \right)^\bell}
    +\displaystyle{\sum_{\bell  \in \NN_0^s, \ |\bell|\le k-1}c_{\bj,\balpha,\bell} \sum_{\bbeta\in \ZZ^s}
     \ra_{\balpha-M\bbeta} \left( M^{-r} \bbeta \right)^\bell}\\
    \\
     &=\displaystyle{\sum_{\bell  \in \NN_0^s, \ |\bell|= k}c_{\bj,\balpha,\bell} \sum_{\bbeta\in \ZZ^s}
     \ra_{\balpha-M\bbeta} \left( M^{-r} \bbeta \right)^\bell}+\displaystyle{\sum_{\bell  \in \NN_0^s, \
     |\bell|\le k-1}c_{\bj,\balpha,\bell} \left( M^{-(r+1)} (\balpha-\btau)\right)^\bell}\\
     \\
      &=\displaystyle{\sum_{\bell  \in \NN_0^s, \ |\bell|= k}c_{\bj,\balpha,\bell} \left(\sum_{\bbeta\in \ZZ^s}
       \ra_{\balpha-M\bbeta} \left( M^{-r} \bbeta \right)^\bell-
      \left(M^{-(r+1)} (\balpha-\btau)\right)^\bell\right)}+q_{\bj}(\btau)\,.
  \end{array}
    $$
    The last equality is due to the fact that
    $$
      q_\bj(\btau)=q_\bj\left(\balpha-M^{r+1}  M^{-(r+1)} (\balpha-\btau) \right)=\sum_{\bell  \in \NN_0^s, \
      |\bell|\le k} c_{\bj,\balpha,\bell} \left(M^{-(r+1)} (\balpha-\btau)\right)^\bell.
    $$
    Note that all rows and columns of $M^{r+1}$, $r \ge 0$, are non-zero vectors, due to
    $\rho(M)>1$ and $det(M^{r+1})=det(M)^{r+1}$. Thus, if all diagonal elements of $M^{r+1}$ are non-zero,
    $c_{\bj,\balpha, \bj} \not =0$. If one of the diagonal elements of $M^{r+1}$ is zero, then there
    exists an $s \times s$ permutation matrix $P$, such that
    for any $\bj \in \NN_0^s$, $|\bj|=k$, the coefficient $c_{\bj,\balpha,P \bj} \not =0$ in
    $$
    \displaystyle{\sum_{\bell  \in \NN_0^s, \ |\bell|= k}c_{\bj,\balpha,\bell} \left(\sum_{\bbeta\in \ZZ^s}
       \ra_{\balpha-  M \bbeta} \left( M^{-r} \bbeta \right)^\bell-
      \left(M^{-(r+1)} (\balpha-\btau)\right)^\bell\right)}=0.$$
    Hence, comparing the coefficients on both sides of the above identity, we note that it is satisfied  for all
    $\bj \in \NN_0^s$, $|\bj|=k$, if and only if
$$ \sum_{\bbeta\in \ZZ^s} \ra_{\balpha-M\bbeta} \left( M^{-r}
\bbeta \right)^\bj-
      \left(M^{-(r+1)} (\balpha-\btau)\right)^\bj=0,\quad \hbox{for all}\  \bj \in \NN_0^s, \quad
      |\bj|=k. \qquad \hbox{\eop}
$$

\smallskip \noindent
The following result follows
directly from \cite[Theorem 2.6]{CharinaConti2012}  and Proposition \ref{prop:equiv3} by replacing
$mI$ by $M\in \ZZ^{s \times s}$ whenever appropriate.

\begin{theorem}\label{teo:main}
Let $k\in \NN_0$. A non-singular  $s$-variate scalar subdivision scheme with symbol
$a(\bz)$ and associated parametrization in (\ref{def:general_par})
with some $\btau \in \RR^s$ reproduces polynomials of degree up to
$k$ if and only if
$$
 \bigl(D^\bj a \bigr)(\bone)= m \ q_{\bj}(\btau) \quad \hbox{and} \quad \bigl(D^\bj a\bigr)(\beps)=0 \quad
\hbox{for} \quad \beps\in \Xi',\quad |\bj| \le k.$$
\end{theorem}

We conclude by observing that, if a scheme is only convergent, then Theorem \ref{teo:main} leads to the
following sufficient conditions for determining the degree of  polynomial reproduction of the scheme.

\begin{corollary}\label{cor:main}
Let $k\in \NN_0$. A convergent $s$-variate scalar subdivision scheme with symbol
$a(\bz)$ and associated parametrization in (\ref{def:general_par})
with some $\btau \in \RR^s$ reproduces polynomials of degree up to
$k$ if
$$
 \bigl(D^\bj a\bigr)(\bone)= m \ q_{\bj}(\btau) \quad \hbox{and} \quad \bigl(D^\bj a\bigr)(\beps)=0 \quad
\hbox{for} \quad \beps\in \Xi',\quad |\bj| \le k\,.$$
\end{corollary}

\section{Applications and examples}
 \label{sec:examples}

\subsection{Box splines}

The results of \cite{CharinaContiJetterZimm2010} imply that any bivariate convergent subdivision scheme
whose symbol satisfies Condition~Z$_3$ is an affine combination of the $3-$directional box spline symbols
from the set
$$
  \{ B_{j,j,h}, \ B_{j,h,j}, \ B_{h,j,j}\ : \
   h=0,1, \quad j=3-h  \}
$$
with
$$
  B_{h,i,j}(\bz)=4 \cdot  \left(\frac{1+z_1}{2} \right)^h \cdot
  \left(\frac{1+z_2}{2} \right)^i \cdot \left(\frac{1+z_1z_2}{2} \right)^j, \quad h,i,j \in \NN_0.
$$
One easily checks that e.g. the affine combination
$$
 a(\bz)=5 \cdot B_{221}(\bz)- B_{212}(\bz) - B_{122}(\bz)-2 \cdot B_{330}(\bz)
$$
satisfies all conditions of Corollary \ref{cor:main} for $k=3$ although none of
the summands of this affine combination does separately. The subdivision scheme associated
with $a(\bz)$ is $C^1$ as its second difference operator is contractive.

\subsection{$\sqrt{3}-$ subdivision}

This example shows how to determine the degree of polynomial reproduction of
the approximating
$\sqrt{3}-$subdivision schemes given in \cite[page 21]{JiangOswald}
 from the corresponding mask symbol
$a(\bz)$ instead of the iterated symbol $a(z_1z_2^{-2},z_1^2z_2^{-1}) \cdot a(\bz)$, as it is done
in \cite{CharinaConti2012}. We also show how to use affine combinations of these schemes to improve their
degree of polynomial reproduction.

The dilation matrix in this case is
$$
  M=\left( \begin{array}{rr} 1&2\\-2&-1\end{array}\right), \quad M^2=-3I,
$$
and, e.g., the mask symbol is given by
\begin{eqnarray}
 a(\bz)&=&\frac{1}{6}\left(z_1z_2+z_1^{-1}z_2^{-1}+z_1^{-1}z_2^2+z_1^{-2}z_2+z_1z_2^{-2}+z_1^2z_2^{-1}\right)
 \notag \\ &+&\frac{1}{3} \left(z_1^{-1}+z_2+z_1z_2^{-1}\right)+\frac{1}{3}
 \left(z_2^{-1}+z_1+z_1^{-1}z_2\right). \notag
\end{eqnarray}
The associated subdivision scheme satisfies zero conditions of at most
order $2$, see \cite{JiangOswald}. The result of Corollary
\ref{cor:linear-reproduction} yields $\btau=(0,0)$, which implies
that the corresponding scheme reproduces linear polynomials, if we
set $\btau=(0,0)$ in (\ref{def:general_par}). Since, the mask
symbol satisfies at most zero conditions of order $2$, the associated
refinable function has approximation order $2$, see \cite{ALevin2003}. Note that,
similarly, the corresponding $\btau$ is $(0,0)$ for all
approximating $\sqrt{3}-$subdivision schemes given in \cite[page
21]{JiangOswald}.

Take next the symbols $a_j(\bz)$, $j=1,2,3,4$ of the four approximating subdivision schemes in \cite[page
21]{JiangOswald} that satisfy zero conditions of order $3$ and consider their affine combination
$$
 a(\bz)=\sum_{j=1}^4 \lambda_j \cdot a_j(\bz), \quad \sum_{j=1}^4 \lambda_j=1.
$$
The resulting scheme still satisfies Condition~Z$_3$ and in addition the rest of the conditions in
Corollary \ref{cor:main} for $k=3$, if
$$
 \lambda_1=\frac{1}{2} \lambda_3+\lambda_4+2, \quad \lambda_2=-\frac{3}{2} \lambda_3-2\lambda_4-1,\quad
 \lambda_3,\lambda_4 \in \RR.
$$
Surprisingly, for any such choice of the parameters $\lambda_1, \dots, \lambda_4$ the resulting
scheme is given by
$$
 a(\bz)=1-\frac{1}{9}(z_1^{-2}+z_1^{-2}z_2^2+z_2^2+z_1^2+z_1^2z_2^{-2}+z_2^{-2})+
 \frac{4}{9}(z_1^{-1}+z_1^{-1}z_2+z_2+z_1+z_1z_2^{-1}+z_2^{-1})
$$
and defines the interpolatory scheme considered in \cite[page 17]{JiangOswald}. Thus, by \cite{ALevin2003},
this scheme has approximation order  $3$.

\subsection{General expanding dilation}

We use the results of \cite{GH} to define a convergent subdivision scheme with dilation matrix
\begin{equation} \label{def:general_M}
   M=\left( \begin{array}{cc} 2&1\\0&2\end{array}\right).
\end{equation}
Similarly, we could also treat the case of any general integer
expanding dilation matrix. As shown in \cite{GH}, one particular
choice for the coset representatives of any such $M$ is obtained by
determining the intersection of $M \cdot [0,1)^2$ with $\ZZ^2$.
Thus, for $M$ in \eqref{def:general_M} we get
$ \rE=\{(0,0),(1,0),(1,1),(1,2)\}$. Applying
the strategy described in \cite{CDM91,Dahlke} we get via
convolution for the mask symbol
$$
 a(\bz)=\sum_{\be \in \rE} \bz^\be
$$
a convergent scheme associated with the symbol $\frac{1}{4}a^2(\bz)$. The symbol $a^2(\bz)$ satisfies
\eqref{eq:Zkcond} with $k=1$ and $\Xi=\{(1,1), (1,-1), (-1,{\mathrm i}), (-1,-{\mathrm i})\}$, ${\mathrm i}=\sqrt{-1}$. The conditions
in \eqref{eq:Zkcond} are not satisfied for $k=2$, e.g. for $\bj=(0,2)$ and $\varepsilon=(1,-1)$.
Using the result of Corollary  \ref{cor:linear-reproduction}
 we get $\btau=(2,1)$, which
together with the zero conditions of order $2$, by \cite{ALevin2003}, implies that the scheme has approximation order $2$.


\end{document}